\documentclass[11pt]{article}
\usepackage{latexsym, amsfonts, amssymb, amsmath, amsthm, graphicx}                              
\def\~{\tilde }
\def\Q{{\sf Q}}
\def\P{{\sf P}}
\def\QR{\overleftarrow{\sf Q}}
\def\PR{\overleftarrow{\sf P}}

\def\psiR{\overleftarrow{\psi}}
\def\PsiR{\overleftarrow{\Psi}}

\def\pp{{\sf p}}

\def\MM{{\cal M}}

\theoremstyle{plain}
\newtheorem{theorem}{Theorem}[section]
\newtheorem{lemma}[theorem]{Lemma}
\newtheorem{corollary}[theorem]{Corollary}

\theoremstyle{definition}
\newtheorem{definition}[theorem]{Definition}
\newtheorem{example}[theorem]{Example}

\begin{document}

\title{Vertex and edge expansion properties for rapid mixing}

\author {Ravi Montenegro \thanks{School of Mathematics, Georgia Institute of Technology, 
         Atlanta, GA 30332-0160, ravi.montenegro@math.gatech.edu; supported in part by
         a VIGRE grant.}}

\date{}

\maketitle
\begin{abstract}
\noindent
We show a strict hierarchy among various edge and vertex expansion properties
of Markov chains. This gives easy proofs of a range of bounds, both classical 
and new, on chi-square distance, spectral gap and mixing time.
The $2$-gradient is then used to give an isoperimetric proof that a random walk on
the grid $[k]^n$ mixes in time $O^*(k^2\,n)$.

\vspace{1.5ex}
\end{abstract}

\section{Introduction}

Markov chain algorithms have been used to solve a variety of previously intractable 
approximation problems. These have included approximating the permanent, estimating volume, counting
contingency tables, and studying stock portfolios, among others. 
In all of these cases a critical point has been to show that a Markov chain is rapidly mixing,
that is, within a number of steps polynomial in the problem size the Markov chain approaches a stationary 
(usually uniform) distribution $\pi$.

Intuitively, a random walk on a graph (i.e., a Markov chain) is rapidly mixing if there are no bottlenecks.
This isoperimetric argument has been formalized by various authors. Jerrum and Sinclair \cite{JS88.1} 
showed rapid mixing occurs if and only if the underlying graph has sufficient edge expansion,
also known as high {\em conductance}. Lov\'asz and Kannan \cite{LK99.1} showed that the
mixing is faster if small sets have larger expansion. Kannan, Lov\'asz and Montenegro \cite{KLM03.1}
and Morris and Peres \cite{MP03.1} extended this and showed that the mixing is even faster 
if every set also has a large number of boundary vertices, i.e., good vertex expansion.

In a separate paper \cite{Mon03.1} the present author has shown that the extensions 
of \cite{LK99.1,KLM03.1,MP03.1} almost always improve on the bounds of \cite{JS88.1},
by showing that standard methods used to study conductance --
via geometry, induction or canonical paths -- can be extended to show that
small sets have higher expansion or that there is high vertex expansion.
This typically leads to bounds on mixing time that are within a single order of magnitude from optimal.
However, none of these methods fully exploit the results of \cite{KLM03.1,MP03.1}, 
as each involves only two of three properties: edge expansion, vertex expansion and conditioning on set size.

Before introducing our results, let us briefly discuss the measures of set expansion / congestion that are 
used in \cite{KLM03.1,MP03.1}. Note that for the remainder of the paper {\em congestion} and {\em bottleneck} 
mean that there are either few edges from a set $A$ to its complement, or that there are few boundary vertices,
i.e. either edge or vertex expansion is poor. 
Kannan et. al. developed {\em blocking conductance} bounds on mixing for three measures of congestion.
The {\em spread} $\psi^+(x)$ measures the worst congestion for sets of sizes in $[x/2,x]$,
so if there are bottlenecks at small set sizes but not at larger ones then this is
a good measure to use. In contrast, the {\em modified spread} $\psi_{mod}(x)$ measures
the worst congestion among sets of all sizes $\leq x$, but comes with stronger mixing bounds, 
so for a ``typical'' case where the congestion gets worse as set
size increases then this is best. The third measure, {\em global spread} $\psi_{gl}(x)$ measures a weighted 
congestion among sets of sizes $\leq x$ which is best only if the Markov chain has
extremely low congestion at small sets. Finally, Morris and Peres' {\em evolving sets}
uses a different measure $\psi_{evo}(x)$ of the worst congestion among sets of sizes $\leq x$. 
Their method comes with very good mixing time bounds in terms of $\psi_{evo}(x)$, 
and because it bounds the stronger chi-square distance it also 
implies bounds on the spectral gap.
However, it is unclear how the size of $\psi_{evo}(x)$ compares with the three congestion measures 
of blocking conductance and hence we do not know which method is best unless all four of the 
congestion measures are computed, a non-trivial task.

We begin by showing that the spread $\psi^+$ lower bounds the evolving sets quantity 
$\psi_{evo}$ of Morris and Peres \cite{MP03.1}. This implies a non-reversible form of $\psi^+$,
as well as lower bounds on the spectral gap and on chi-square distance. 
The other forms of blocking conductance are found to upper bound $\psi_{evo}$ and are more
appropriate for total variation distance.
Moreover, an ``optimistic'' form of the spread turns out to upper bound the spectral gap
and lower bound total variation mixing time, although this form is not useful in practice.

Houdr\'e and Tetali \cite{HT03.1}, in the context of concentration inequalities, considered  
the discrete gradients $h_p^+(x)$, a family which
involves all three properties of the new mixing methods -- edges, vertices and set size --
with $p=1$ measuring only edges, $p=2$ weighting edges and vertices roughly equally,
and $p=\infty$ measuring only vertices. In this paper it is shown that the spread function 
$\psi^+(x)$ is closely bounded both above and below by $h_p^+(x)$.
It is found that various classical isoperimetric bounds on mixing time and spectral gap
are essentially the best lower bound approximations to the quantity $h_2^+(x)^2/2$. 
The $h_1^+(1/2)$ approximation is the theorem of Jerrum and Sinclair, 
$h_1^+(x)$ leads to the average conductance of Lov\'asz and Kannan,
$h_{\infty}^+(1/2)$ gives a mixing time bound of Alon \cite{Alon86.1}, and $h_2^+(x)$ gives a
bound shown by Morris and Peres \cite{MP03.1} and in a weaker form by this author \cite{Mon02.1}.

Of these various bounds the one that is the most relevant to our purposes is $h_2^+(x)$, since
this is weighted equally between edge and vertex isoperimetry.
In order to give an application of our methods we show how two additional isoperimetric quantities,
Bobkov's constant $b_p^+$ and Murali's $\beta^+$ \cite{Mur01.1}, can be used to bound $h_2^+(x)$ for products
of Markov chains. We apply this to prove a lower bound on $h_2^+(x)$ for a 
random walk on the grid $[k]^n$. This leads to a mixing time bound of $O(k^2\,n\,\log^2 n)$,
the first isoperimetric proof of the correct $\tau = O^*(k^2\,n)$ for this Markov chain.

The paper proceeds as follows. In Section \ref{sec:prelim} we introduce notation.
Section \ref{sec:chi2} shows the connection between spread, evolving sets and spectral gap.
Section \ref{sec:gradients} gives results on the discrete gradients, including
sharpness. Section \ref{sec:applications} finishes
the paper with the isoperimetric bound on the grid $[k]^n$.


\section{Preliminaries}  \label{sec:prelim}

A {\em finite state Markov chain} $\MM$ is given by a state space $K$ with 
cardinality $|K| = n$, and the {\em transition probability matrix}, an 
$n\times n$ square matrix $\P$ such that $\P_{ij} \in [0,1]$ and
\mbox{$\forall i \in K\ :\ \sum_{j \in K} \P_{ij} = 1$}. Probability distributions
on $K$ are represented by $1\times n$ row vectors, so that if the initial 
distribution is $\pp^{(0)}$ then the $t$-step distribution is
given by $\pp^{(t)} = \pp^{(0)}\,\P^t$. 

The Markov chains considered here are {\em irreducible} ($\forall i,j \in K,\ \exists t : (\P^t)_{ij} > 0$)
and {\em aperiodic} ($\forall i : gcd\{ t : (\P^t)_{ii} > 0\} = 1$). Under these conditions there is
a unique {\em stationary distribution} $\pi$ such that $\pi\,\P = \pi$,
and moreover the Markov chain is {\em ergodic} ($\forall i,\,j\in K : \lim_{t\rightarrow \infty} (\P^t)_{ij} = \pi_j$).
All Markov chains in this paper are {\em lazy} ($\forall i\in K:\P_{ii} \geq \frac 12$); 
lazy chains are obviously aperiodic.

The {\em time reversal} of a Markov chain $\MM$ is the Markov chain with transition probabilities
$\PR(u,v) = \pi(v)\P(v,u)/\pi(u)$, and has the same stationary distribution $\pi$ as the original Markov chain.
It is often easier to consider {\em time reversible} Markov chains 
($\forall i,j \in K : \pi(i)\,\P(i,j) = \pi(j)\,\P(j,i)$). In the time reversible case $\PR=\P$ and the reversal
is just the original Markov chain.

The distance of $\pp^{(t)}$ from $\pi$ is measured by the $L^p$ distance $\|\cdot\|_{L^p(\pi)}$,
which for $p\geq 1$ is given by 
$$
\| \pp^{(t)} - \pi \|_{L^p(\pi)}^p = \sum_{v\in K} \left| \frac{\pp^{(t)}(v)}{\pi(v)}-1 \right|^p\,\pi(v) \ .
$$
The {\em total variation distance} is $\|\cdot\|_{TV}=\frac 12\,\|\cdot\|_{L^1(\pi)}$, and the
{\em $\chi^2$-distance} is $\|\cdot\|_{\chi^2(\pi)}=\|\cdot\|_{L^2(\pi)}^2$. 

The {\em mixing time} measures how many steps it takes a Markov chain to approach the stationary distribution,
\begin{eqnarray*}
\tau(\epsilon)
&=& \max_{\pp^{(0)}} \min \left\{ t : \| \pp^{(t)} - \pi \|_{TV} \leq \epsilon \right\}, \\
\chi^2(\epsilon) &=&
  \max_{\pp^{(0)}} \min \left\{ t : \| \pp^{(t)} - \pi \|_{\chi^2(\pi)} \leq \epsilon \right\}\ . 
\end{eqnarray*}
Cauchy-Schwartz shows that $2\,\|\cdot\|_{TV} \leq \|\cdot\|_{\chi^2(\pi)}^{1/2}$, from which it
follows that $\tau(\epsilon) \leq \chi^2(4\epsilon^2)$.

Morris and Peres showed a nice fact about general (non-reversible) Markov chains \cite{MP03.1}
$$
\max_{x,z} \frac{\P^{n+m}_{xz}}{\pi(z)}-1 \leq 
   \|\P^n(x,\cdot)-\pi\|_{\chi^2(\pi)}^{1/2}\,\|\PR^m(z,\cdot)-\pi\|_{\chi^2(\pi)}^{1/2}\, ,
$$
and so chi-square mixing can be used to show small
relative pointwise distance $\pp^{(t)}(\cdot)/\pi(\cdot)$.
This makes chi-square mixing a stronger condition than total variation mixing.

The {\em ergodic flow} between two points $i,j\in K$ is $q(i,j) = \pi_i\,\P_{ij}$ and the flow between 
two sets $A,\,C\subset K$ is $\Q(A,C) = \sum_{\substack{i\in A \\ j\in C}} q(i,j)$. 
In fact $\Q(A,A^c) = \Q(A^c,A)$, where $A^c:=K\setminus A$.

The {\em continuization} $\tilde K$ of $K$ is defined as follows. Let $\tilde K=[0,1]$ and
to each point $v\in K$ assign an interval $I_v=[a,b]\subset[0,1]=\~K$ with $b-a=\pi(v)$,
so that $m(I_x\cap I_y)=0$ if $x\neq y$, and 
$[0,1]$ is the union of these intervals.
Then if $A,\,B\subset \tilde K$ define $\pi(A) = m(A)$ and
$\Q(A,B) = \sum_{x,y\in K} \frac{m(A\cap I_x)}{\pi(x)}\,\frac{m(B\cap I_y)}{\pi(y)}\,q(x,y)$
for Lebesgue measure $m$. This is consistent with the definition of
ergodic flow between sets in $K$. The continuization $\~K$ is somewhat awkward but will
be needed in our work, particularly for $\psi^+(A)$ below and in the next few sections.

Various isoperimetric quantities have been used to upper bound $\tau(\epsilon)$ and $\chi^2(\epsilon)$.
A few of them are listed below. Unless explicitly stated, all sets both here and later in the paper
will be in $K$, not $\~K$.
$$
\begin{array}{llll}
\vspace{1ex}
Conductance \cite{JS88.1,LK99.1} &:&
\multicolumn{2}{l}{
 \displaystyle \Phi(x) = \min_{\pi(A)\leq x} \Phi(A)\ ,
 \quad \displaystyle \Phi(A) = \frac {\Q(A,A^c)}{\pi(A)}\ ,
 \quad \Phi = \Phi(1/2) ,
} \\
\vspace{1ex}
Spread \cite{KLM03.1} &:&
\displaystyle h^+(x)=\sup_{\substack{A\subset\~K,\\x/2\leq \pi(A)\leq x}} \frac{1}{x\,\psi^+(A)} \ &where\  
\displaystyle \psi^+(A)=\int_0^{\pi(A)} \frac{\Psi(t,A^c)}{\pi(A)^2}\,dt , \\
\vspace{1ex}
Modified\ Spread \cite{KLM03.1} &:&
\displaystyle h_{mod}(x)=\sup_{\substack{A\subset\~K,\\ \pi(A)\leq x}} \frac{1}{x\,\psi_{mod}(A)} &where\  
\displaystyle \psi_{mod}(A)=\int_0^1 \frac{\Psi(t,A^c)}{\pi(A)\,t^{\#}}\,dt , \\
\vspace{1ex}
Global\ Spread \cite{KLM03.1} &:&
\displaystyle h_{gl}(x)=\sup_{\substack{A\subset\~K,\\ \pi(A)\leq x}} \frac{1}{\pi(A)\,\psi_{gl}(A)} &where\  
\displaystyle \psi_{gl}(A)=\int_0^1 \frac{\Psi(t,A^c)}{\pi(A)^2}\,dt , \\
\vspace{1ex}
Evolving\ Sets \cite{MP03.1} &:&
\displaystyle \psi_{evo}(x)=\inf_{\pi(A)\leq x} \psi_{evo}(A) &where\ 
\displaystyle \psi_{evo}(A)=1-\int_0^1 \sqrt{\frac{\pi(A_u)}{\pi(A)}}\,du , 
\end{array}
$$
$$
\begin{array}{l}
\displaystyle \qquad where\ 
A_u=\{y\in K : \PR(y,A) > u\} 
\\
\qquad\ and\ 
\Psi(t,A^c) =
\begin{cases}
\inf_{\substack{\~K\supset B\subset A,\\ \pi(B)= t}} \Q(B,A^c) & if \ t\leq \pi(A)\\
\Psi(1-t,A) & if\ t>\pi(A)
\end{cases}
\hspace{1in}
\end{array}
$$
and $t^{\#}=\min\{t,1-t\}$.

Properties of the various $\psi$ quantities can be found in the introduction and in the following section.

The infinum in the definition of $\Psi(t,A^c)$ is achieved for some set $B\subset \tilde K$ when the 
Markov chain is finite. For instance, when $t\leq\pi(A)$ then one construction for $B$ is 
as follows: order the points $v\in A$ in increasing order of $\P(v,A^c)$ as $v_1,\,v_2,\,\ldots,\,v_n$, 
add to $B$ the longest initial segment $B_m:=\{v_1,\,v_2,\,\ldots,\,v_m\}$ of these points with 
$\pi(B_m)\leq t$, and for the remainder of $B$ take $t-\pi(B_m)<\pi(v_{m+1})$ units of $v_{m+1}$.

The quantities $A_u$ and $\Psi(t,A^c)$ are closely related. For time-reversible
Markov chains, if $t=\pi(A_u)\leq\pi(A)$ 
then $A_u$ is the set of size $t$ with the highest flow into $A$, so
the smallest flow into $A^c$, and therefore $\Psi(t,A^c)=\Q(A_u,A^c)$.
Similarly, when $t=\pi(A_u)>\pi(A)$ then $\Psi(t,A^c)=\Q(A_u^c,A)$.
The choice of $\pi(A_u)$ or $\Psi(t,A^c)$ is similar to the choice of Lebesgue
or Riemann integral, where the Lebesgue-like $\pi(A_u)$ measures the amount of $K$
with transition probability to $A$ above a certain level, while $\Psi(t,A^c)$ is more Riemann-like in
simply integrating along $\P$ once $\P$ has been put in increasing order.

The quantities $\Phi$ and $\psi_{evo}(A)$ can be used to upper bound
$\chi^2(\epsilon)$, while $\Phi(A)$, $\psi^+(A)$ and $\psi_{gl/mod}(A)$ are
used to upper bound $\tau(\epsilon)$. In the $\tau(\epsilon)$ case it suffices to
bound $\tau(1/4)$, because $\tau(\epsilon)\leq \tau(1/4)\,\log_2(1/\epsilon)$ \cite{AF.1}.
The two bounds of most interest here are:

\begin{theorem} \label{thm:mixingtime}
If $\MM$ is a (lazy, aperiodic, ergodic) Markov chain then
\begin{eqnarray*}
\cite{KLM03.1} \quad& \tau(1/4)  
    &\leq 8\cdot1376\,\left(\int_{\pi_0}^{1/2} h(x)\,dx+h(1/2)\right) \\
\cite{MP03.1} \quad&
\chi^2(\epsilon) &\leq \int_{\pi_0}^{1/2} \frac{dx}{x\,\psi_{evo}(x)} + \frac{\log(8/\epsilon)}{\psi_{evo}(1/2)}
\end{eqnarray*}
where $\pi_0=\min_{v\in K} \pi(v)$ and $h(x)$ indicates $h^+(x)$, $h_{mod}(x)$
or $h_{gl}(x)$. The $h(x)$ bounds apply to reversible Markov chains only, whereas the
$\psi_{evo}(x)$ bound applies even in the non-reversible case.
\end{theorem}


\section{Spread, $\chi^2$ and the Spectral Gap} \label{sec:chi2}

In this section we show a connection between the spread function and evolving sets.
We further explore this connection by finding that variations on the spread function both
upper and lower bound the spectral gap. The connection to evolving sets
implies a mixing time theorem with much stronger constants, as well as
a non-reversible result.

\begin{theorem} \label{thm:main}
If $\MM$ is a lazy Markov chain and $A$ is a subset of the state space with $\pi(A)\leq 1/2$,
then let the (time reversed) spread function $\psiR^+(A)$ be given by
$$
\psiR^+(A) = \int_0^{\pi(A)} \frac{\PsiR(t,A^c)}{\pi(A)^2}\,dt
\quad where \quad
\PsiR(t,A^c) = \inf_{\substack{\~K\supset B\subset A,\\ \pi(B)=t}} \Q(A^c,B) \ ,
$$
with $\psiR_{gl}(A)$ and $\psiR_{mod}(A)$ defined similarly, and 
where $\psiR^+(x)=\inf_{\pi(A)\leq x} \psiR^+(A)$. Then
$$
\psiR_{gl}(A) \geq \frac 12\,\psiR_{mod}(A) \geq \psi_{evo}(A) \geq \frac 14\,\psiR^+(A)\ ,
$$
and in particular, 
$$
\tau(\epsilon) \leq \chi^2(4\,\epsilon^2) \leq
\begin{cases}
\vspace{2ex}
\displaystyle 4\,\int_{\pi_0}^{1/2} \frac{dx}{x\,\psiR^+(x)} + \frac{4\,\log(2/\epsilon^2)}{\psiR^+(1/2)} & \\
\displaystyle \frac{2}{\psiR^+(1/2)}\,(\log (1/\pi_0) + 2\,\log(1/2\epsilon))\ .
\end{cases}
$$
\end{theorem}

Observe that $\psiR^+$ is just $\psi^+$ of the time reversal. In particular,
$\psiR^+(A)=\psi^+(A)$ when the Markov chain is reversible,
so this is an extension of the results of Kannan et. al. \cite{KLM03.1}.

Corollary \ref{thm:extensions} shows that $\psiR^+(1/2)\geq \Phi^2$ for
lazy Markov chains, because $\overleftarrow{\Phi}=\Phi$ even for non-reversible
Markov chains. This approximation applied to the second upper bound on $\tau(\epsilon)$ 
is exactly a factor two from the non-reversible bound shown by Mihail \cite{Mih89.1}.
A more direct approach can be found in \cite{MP03.1} which recovers this factor of two.

The inequalities
$\psiR_{gl}(A)\geq \psiR^+(A)$ and $\psiR_{mod}(A)\geq \psiR^+(A)$
follow almost immediately from the definitions, so the theorem should
not be used to lower bound either of these quantities by $\psiR^+(A)$.

Nevertheless, the inequalities between $\psi$ terms given in the theorem are all sharp. 
The first two inequalities are sharp for a walk with uniform transition probability 
$\alpha/2\leq 1$ from $A$ and all the flow concentrated in a region of size 
$\alpha\pi(A)$ in $A^c$. The final inequality is sharp as a limit.
Let $D\rightarrow 4^-$, $x_0=(D/4)^{2/3}$, $\alpha=(4/D)-1$.
Then put an $1-x_0$ fraction of $A$ with $\P(\cdot,A^c)=\alpha/2$
and the remainder with $\P(\cdot,A^c)=0$. This flow can be concentrated 
in a small region of $A^c$.

Even though $\psi_{gl}(A)$ is the largest quantity it is usually the least useful. 
As discussed in the introduction, when there are bottlenecks at small values of 
$\pi(A)$ then $h^+(x)$ is best (i.e., smallest) because of the conditioning on $\pi(A)\in[x/2,x]$.
Spread $\psi^+(A)$ is also the easiest to compute, and the connection to $\psi_{evo}(A)$
improves the constant terms in Theorem \ref{thm:mixingtime} greatly. For a ``typical'' case 
$h_{mod}(x)$ is better than $h^+(x)$, but $h_{gl}(x)$ is poor because the
supremum in $h_{gl}(x)$ may occur for small $\pi(A)$. However, for graphs with 
extremely high node-expansion then $h_{gl}(x)$ may be best.
As a case in point, on the complete graph $K_n$ we have $\tau(1/4)\leq \chi^2(1/4) = O(\log n)$
via $\psi^+$ or $\psi_{evo}$, while $\tau(1/4)=O(\log \log n)$ from $\psi_{mod}$ and
$\tau(1/4)=O(1)$ from $\psi_{gl}$. However, on the cube $\{0,1\}^n$,
$\psi_{gl}$ implies only $\tau=O(n2^n)$, hopelessly far from the correct
$\tau=O(n\log n)$.

The following lemma shows how to rewrite $\PsiR(t,A^c)$ in terms of $\pi(A_u)$ and will be key to our proof.

\begin{lemma} \label{lem:main}
If $\MM$ is a lazy Markov chain and $A\subset K$ is a subset of the state space then
$$
\PsiR(t,A^c) = 
\begin{cases}
\vspace{1ex}
\displaystyle \int_{w(t)}^1 (t-\pi(A_u))\,du & if\ t\leq \pi(A) \\
\displaystyle \int_0^{w(t)} (\pi(A_u)-t)\,du & if\ t> \pi(A) 
\end{cases}
\quad where \quad
w(t) = \inf\{y\,:\,\pi(A_y)\leq t\}.
$$
\end{lemma}

\begin{proof}
We consider the case of $t\leq\pi(A)$. A similar argument implies the result when $t>\pi(A)$.

When $t=\pi(A_x)$ for some $x\in[1/2,1]$ then $A_x$ is the set of size $t$ with the highest flow from $A$,
so the smallest flow from $A^c$, and therefore $\PsiR(t,A^c)=\Q(A^c,A_x)$
($\PsiR$ considers the reversed chain, so it minimizes flow from $A^c$ rather than into $A^c$).
But 
$$
\begin{array}{rclcl}
\vspace{1ex}
\displaystyle \Q(A^c,A_x)
  &=& \displaystyle \sum_{y\in A_x} \pi(y)\,\PR(y,A^c) 
  &=& \displaystyle \sum_{y\in A_x} \int_0^1 \left(1-1_{\PR(y,A)\geq u}\right)\,du\,\pi(y) \\
\vspace{1ex}
  &=& \displaystyle \int_0^1 \sum_{y\in A_x} \left(1-1_{\PR(y,A)\geq u}\right)\,\pi(y)\,du
  &=& \displaystyle \int_0^1 \left(\pi(A_x)-\pi(A_u\cap A_x)\right)\,du \\
  &=& \displaystyle \int_{w(t)}^1 \left(\pi(A_x)-\pi(A_u)\right)\,du
  &=& \displaystyle \int_{w(t)}^1 (t-\pi(A_u))\,du\,.
\end{array}
$$
This gives the result if $t = \pi(A_x)$ for some $x$. Otherwise, the set $B$
where the infinum is achieved in the definition of $\PsiR(t,A^c)$ contains
$A_{w(t)+\delta}$ where $\delta\rightarrow 0^+$, and the remaining points
$y\in B\setminus A_{w(t)+\delta}$ satisfy $\PR(y,A)=w(t)$. Let
$x=w(t)+\delta$, then $w(\pi(A_x))=w(t)$ for $\delta$ sufficiently
small and
\begin{eqnarray*}
\PsiR(t,A^c) &=& \QR(A_x,A^c)+\QR(B\setminus A_x,A^c) \\
   &=& \int_{w(\pi(A_x))}^1 \left(\pi(A_x) - \pi(A_u)\right)\,du 
       + (t-\pi(A_x))\,(1-w(t)) \\
   &=& \int_{w(t)}^1 \left(t - \pi(A_u)\right)\,du\,.
\end{eqnarray*}
\end{proof}

\begin{proof}[Proof of theorem]
Rewriting $\psiR^+(A)$ in terms of $\pi(A_u)$ gives
$$
\begin{array}{rclcl}
\vspace{1ex}
\displaystyle \psiR^+(A) 
  &=& \displaystyle \int_0^{\pi(A)} \frac{\PsiR(t,A^c)}{\pi(A)^2}\,dt
  &=& \displaystyle \int_0^{\pi(A)} \int_{w(t)}^1 \frac{t-\pi(A_u)}{\pi(A)^2}\,du\,dt \\
  &=& \displaystyle \int_{1/2}^1 \int_{\pi(A_u)}^{\pi(A)} \frac{t-\pi(A_u)}{\pi(A)^2}\,dt\,du  
  &=& \displaystyle  \frac 12\,\int_{1/2}^1 \left(\frac{\pi(A)-\pi(A_u)}{\pi(A)}\right)^2\,du \ ,
\end{array}
$$
where the first equality follows from the definition of $\psiR^+(A)$, the second
equality applies Lemma \ref{lem:main}, the third is a change in the order of
integration using that $w(t)\leq u$ iff $\pi(A_u) \leq t$, and the final equality
is integration with respect to $t$. 
Morris and Peres \cite{MP03.1} used a Taylor approximation 
$\sqrt{\pi(A_u)/\pi(A)}=\sqrt{1+x}\leq 1+x/2-(x^2/8)\,\delta_{x\leq 0}$
for $x=\pi(A_u)/\pi(A)-1$, and the Martingale property of $\pi(A_u)$ that 
$\int_0^1 \pi(A_u)\,du=\pi(A)$ (Lemma 6 of \cite{MP03.1}), to derive the lower bound 
$$
\psi_{evo}(A) \geq \frac 18\,\int_{1/2}^1 \left(\frac{\pi(A)-\pi(A_u)}{\pi(A)}\right)^2\,du .
$$
The lower bound $\psi_{evo}(A)\geq \psi^+(A)/4$ follows.

Similarly, 
\begin{eqnarray*}
\psiR_{mod}(A) 
  &\geq& \int_0^1 \frac{\PsiR(t,A^c)}{t\,\pi(A)}\,dt \\
  &=& \int_0^{\pi(A)} \int_{w(t)}^1 \frac{t-\pi(A_u)}{t\,\pi(A)}\,du\,dt 
     + \int_{\pi(A)}^1 \int_0^{w(t)} \frac{\pi(A_u)-t}{t\,\pi(A)}\,du\,dt \\ 
  &=& \int_{1/2}^1  \int_{\pi(A_u)}^{\pi(A)} \frac{t-\pi(A_u)}{t\,\pi(A)}\,dt\,du  
     + \int_0^{1/2} \int_{\pi(A)}^{\pi(A_u)} \frac{\pi(A_u)-t}{t\,\pi(A)}\,dt\,du \\
  &=& \int_0^1 \frac{\pi(A_u)}{\pi(A)}\log\left(\frac{\pi(A_u)}{\pi(A)}\right)\,du.
\end{eqnarray*}
The final equality used the Martingale property of $\pi(A_u)$, as does the equality below.
\begin{eqnarray*}
\psi_{evo}(A) &=& \int_0^1 \left(\frac{\pi(A_u)}{\pi(A)} - \sqrt{\frac{\pi(A_u)}{\pi(A)}}\right)\,du \\
 &\leq& \frac 12\,\int_0^1 \frac{\pi(A_u)}{\pi(A)}\log\left(\frac{\pi(A_u)}{\pi(A)}\right)\,du
   \leq \psiR_{mod}(A)/2,
\end{eqnarray*}
where the inequality follows from $2(x-\sqrt{x}) \leq x\log x$ for all $x> 0$.

To establish that $2\psi_{gl}(A)\geq \psi_{mod}(A)$, observe that when $t\in [\pi(A),1-\pi(A)]$ then the
result is trivial, as $1/\min\{t,1-t\}\pi(A) \leq 1/\pi(A)^2$. When $t\in [0,\pi(A)]$ then let $f(t)=\PsiR(t,A^c)/t$.
If $B$ is the set where the infinum in $\PsiR(t,A^c)$ is achieved then $f(t)$
is the average probability of the reversed chain making a transition from a point in $B$ 
to $A^c$ in a single
step. It follows that $f(t)$ is an increasing function, because as $t$ increases
the points added to $B$ will have higher probability of leaving then any
of those previously added. We then have the following:
\begin{eqnarray*}
\int_0^{\pi(A)} \left[\frac{2\,\PsiR(t,A^c)}{\pi(A)^2}-\frac{\PsiR(t,A^c)}{t\,\pi(A)}\right]\,dt
  &=& \int_0^{\pi(A)} \frac{(2t-\pi(A))\,f(t)}{\pi(A)^2}\,dt  \\
 &=& \int_{\pi(A)/2}^{\pi(A)} \frac{(2t-\pi(A))\,(f(t)-f(\pi(A)-t))}{\pi(A)^2}\,dt \geq 0.
\end{eqnarray*} 
A similar argument holds for the interval $t\in[1-\pi(A),1]$.

The first upper bound for $\tau$ follow from 
$4\,\|\cdot\|_{TV}^2\leq\|\cdot\|_{\chi^2(\pi)}$ and Theorem \ref{thm:mixingtime}.
The second follows from this and
$\chi^2(\epsilon)\leq (2\psi_{evo}(1/2))^{-1}\,\log(1/\epsilon\pi_0)$, which
is another bound of Morris and Peres \cite{MP03.1}.
\end{proof}

The connection between the spread function and mixing quantities is deeper than just an
upper bound on mixing time.
In the proof that $\psi^+$ bounds mixing time \cite{KLM03.1} it is shown that for
reversible Markov chains there is
some ordering of points in the state space $\tilde K=[0,1]$
such that the mixing time is lower bounded by the case when 
$\psi_{correct}(x) = \psi_{correct}([0,x]) = \int_0^x \Q([x-t,x],[x,1])\,dt$.
The most pessimistic lower bound on $\psi_{correct}(x)$ is $\psi^+(x)$, hence
an upper bound on mixing time, whereas 
\begin{equation} \label{eqn:big}
\psi_{big}(A)= \int_0^{\pi(A)} \frac{\Psi_{big}(t,A^c)}{\pi(A)^2}\,dt
\quad where\quad \Psi_{big}(t,A^c) = \sup_{\substack{\~K\supset B\subset A,\\ \pi(B)=t}} \Q(B,A^c)
\end{equation}
is the most pessimistic upper bound on $\psi_{correct}(A)$ when $A=[0,x]$, i.e.,
$\psi_{big}(x)\geq \psi_{correct}(x)\geq \psi^+(x)$.
The following theorem shows that this ordering carries over to mixing time and spectral gap,
with $\psi_{big}$ appearing in a lower bound on mixing time and in an upper bound on
spectral gap.

\begin{theorem} \label{thm:spectral}
If $\MM$ is a lazy, aperiodic, ergodic reversible Markov chain then
\begin{eqnarray*}
\frac{1-4\psi_{big}(1/2)}{8\,\psi_{big}(1/2)}\,\log(1/2\epsilon)
  &\leq \tau(\epsilon) \leq& \frac{2}{\psi^+(1/2)}\,(\log (1/\pi_0) + 2\,\log(1/2\epsilon)) \\
4\,\psi_{big}(1/2) &\geq \lambda \geq& \frac 14\,\psi^+(1/2)\ .
\end{eqnarray*}
\end{theorem}

The lower bound on $\lambda$, when combined with $\psi^+(1/2)\geq\Phi^2$ is a factor
two from the well known $\lambda\geq\Phi^2/2$ \cite{Sin92.1}.

\begin{proof}
The upper bound for $\tau(\epsilon)$ follows from the previous theorem.

For the lower bound on $\lambda$ we need some information about the proofs
of certain useful facts. First,  
\begin{equation}\label{eqn:Sin}
\tau(\epsilon) \geq \frac 12\,(1-\lambda)\,\lambda^{-1}\,\ln(2\epsilon)^{-1}\quad(\textrm{see \cite{Sin92.1}})
\end{equation}
can be proven by first showing that $c\,(1-\lambda)^t \leq \sup_{\pp^{(0)}} \|\pp^{(t)}-\pi\|_{TV}$
for some constant $c$. Second, the proof of the $\psi_{evo}$ part of Theorem \ref{thm:mixingtime}
can be easily modified to show that $\sup_{\pp^{(0)}} \|\pp^{(t)}-\pi\|_{TV} \leq c_2\,(1-\psi_{evo}(1/2))^t$
for some constant $c_2$. It follows that
$c\,(1-\lambda)^t \leq c_2\,(1-\psi_{evo}(1/2))^t$, and taking $t\rightarrow\infty$ implies further
that $1-\lambda\leq 1-\psi_{evo}(1/2)$. The result then follows by $\psi_{evo}(A)\geq\frac 14\,\psi^+(A)$.
In words, this says that the asymptotic rate of convergence of total variation
distance is at best $1-\lambda$ and at worst $1-\psi_{evo}(1/2)$, and therefore
$1-\lambda\leq 1-\psi_{evo}(1/2)$.

For the upper bound on $\lambda$ suppose that $\~K\supset A=[0,x]$ where $x=\pi(A)\leq 1/2$, 
and order vertices in $A$ by increasing $\P(t,A^c)$. Then
\begin{eqnarray*}
\psi_{big}(A) + \psi^+(A)
&=& \int_0^{\pi(A)} \frac{\Q([x-t,x],[x,1])}{\pi(A)^2}\,dt 
  + \int_0^{\pi(A)} \frac{\Q([0,t],[x,1])}{\pi(A)^2}\,dt \\
&=& \frac{\Q(A,A^c)}{\pi(A)} = \Phi(A)\ .
\end{eqnarray*}
It follows that $\Phi(A) \geq \psi_{big}(A) \geq \Phi(A)/2$.
The upper bound on $\lambda$ then follow from $\lambda\leq 2\Phi$ \cite{Sin92.1}.

The lower bound on mixing time follows from the upper bound on the
spectral gap and the lower bound on mixing time given in (\ref{eqn:Sin}).
\end{proof}

It would be interesting to know if the lower bound on the mixing time
can be improved. The barbell consisting of two copies of $K_n$ joined by
a single edge is a case where $\tau(1/4) < 1/\psi^+(1/2)$, which
shows that $\psi^+$ cannot replace $\psi_{big}$ in the lower bound. 
However, in those examples where we know the answer 
we find that $\tau(1/4) \geq c\,\int dx/\psi^+(x)$.


\section{Discrete Gradients} \label{sec:gradients}

In this section we look at the discrete gradients $h_p^{\pm}(A)$
of Houdr\'e and Tetali \cite{HT03.1}. This is a family that extends the ideas of
edge and vertex-expansion, with $h_1^{\pm}(A)$ measuring edge-expansion, 
$h_{\infty}^{\pm}(A)$ measuring vertex-expansion and $h_2^{\pm}(A)$ a hybrid. 
We use the $h_p^{\pm}$ notation here, despite the similarity to the
$h_{gl/mod}$ notation earlier, to be consistent with \cite{HT03.1,KLM03.1}.

\begin{definition}
Let $\MM$ be a Markov chain. Then for $p\geq 1$, $A\subset K$ the {\em discrete $p$-gradient $h_p^+(A)$} is
\begin{equation*}
h_p^+(A) = \frac{\Q_p(A,A^c)}{\min\{\pi(A),\pi(A^c)\}} \quad 
where \quad \Q_p(C,D) = \sum_{v\in C} \sqrt[p]{\P(v,D)}\,\pi(v)\ . 
\end{equation*}
The (often larger) $h_p^-(A)$ is defined similarly, but with
$\Q_p(A^c,A)$ rather than $\Q_p(A,A^c)$.

These can be extended to $p=\infty$ in the natural way, by taking
$\Q_{\infty}(C,D) = \pi(\{u\in C : \Q(u,D)\neq 0\})$. 
We sometimes refer to $h_p^{\pm}(x)=\inf_{\pi(A)\leq x} h_p^{\pm}(A)$.
\end{definition}

The main focus of this section will be $h_2^+(A)$ and
$h_2^-(A)$ which are hybrids of edge and vertex-expansion as Cauchy-Schwartz shows,
$h_2^+(A)^2 \leq h_{\infty}^+(A)\,h_1^+(A)$.
Note that $h_2^-(A)$ can be significantly larger than $h_2^+(A)$, which is 
why our theorems below differ in 
the plus and minus cases. In contrast, conductance bounds only 
have one form, for $\Phi(A)=h_1^+(A)=h_1^-(A)$.

In this section it will be shown how discrete gradients can be used to upper and lower bound the chain of 
inequalities given in Theorem \ref{thm:main}. It is not necessary to prove a theorem for the time-reversal
because, for instance, bounds on $\psi^+$ apply to $\psiR^+$ as well by bounding $\psi^+$ of the
time-reversed Markov chain. If we let $\psi^-(A)$ be defined as 
$$
\psi^-(A) = \int_{\pi(A)}^1 \frac{\Psi(t,A^c)}{\pi(A)^2}\, 
$$
then $\psi_{gl}=\psi^+(A) + \psi^-(A)$, so upper and lower bounds on $\psi^{\pm}(A)$ imply
upper and lower bounds on $\psi_{gl}(A)$ as well. Our main result of this section is the following theorem.

\begin{theorem} \label{thm:bounds}
Given a (non-reversible) Markov chain $\MM$ with state space $K$ and a set $A\subset K$, let 
$\P_* = 1 - \inf_{u\in A} \P(u,A)$ and 
$\displaystyle \P_{min} = \inf_{\substack{u\in A,v\in A^c,\\ \P(u,v)>0}} \frac{\P(u,v)}{\pi(v)}$. 
If $\pi(A)\leq \frac 12$ then
\vspace{1ex}

$$
\begin{array}{rcccl}
\vspace{1ex}
\displaystyle \frac 12\,h_2^{\pm}(A)^2 &\geq&\psi^{\pm}(A)&\geq&
    \displaystyle \frac 12\,h_2^{\pm}(A)^2\left/\log\left(\frac{12\,h_1^{\pm}(A)\,h_{\infty}^{\pm}(A)}{h_2^{\pm}(A)^2}\right)\right.  \\
\vspace{1ex}&& and &\geq& 
    \displaystyle \frac 12\,h_2^{\pm}(A)^2\,\sqrt{\frac{\P_{min}}{\P_*}}\ . \\
\end{array}
$$
\end{theorem}

In practice it may be useful to upper bound the $\log$ terms either
by $\log (12\,\P_*/\P_{min})$ for $\psi^{\pm}$, or with 
$\log (12/h_2^+(A)^2)$ for the $\psi^+(A)$ case and 
$\log (12/\pi(A)\,h_2^-(A)^2)$ for $\psi^-(A)$. These follow by the
identities $h_1^{\pm}(A)\leq h_{\infty}^{\pm}(A)\,\P_*$ and 
$h_2^{\pm}(A)\geq h_{\infty}^{\pm}(A)\,\sqrt{\P_{min}}$ for the first type, or by
$h_1^{\pm}(A),\,h_{\infty}^+(A) \leq 1$ and $h_{\infty}^-(A)\leq \pi(A)^{-1}$
for the latter. 

All upper and lower bounds scale properly in $\P$, e.g., if $\P$ is slowed by a factor of $2$ to
$\P \rightarrow \frac 12\,(I+\P)$ then $\psi^{\pm}(A)$ and all the bounds in 
Theorem \ref{thm:bounds} change by the same factor of $2$. Moreover, if $\P(\cdot,A^c)$ is constant over a set $A$ 
then the upper bounds are sharp, while the lower bounds are within a 
small constant factor.

Our methods also extend to the other discrete gradients $h_p^{\pm}$. 
The most interesting cases are $p=1$ and $p=\infty$.

\begin{theorem} \label{thm:extensions}
Given the conditions of Theorem \ref{thm:bounds} then
$$
\frac 12\,h_1^{\pm}(A)\,h_{\infty}^{\pm}(A) \geq \psi^{\pm}(A) \geq
\max\left\{
         \frac 12\, \P_{min}\,h_{\infty}^{\pm}(A)^2,\ 
         \frac 12\,\frac{h_1^{\pm}(A)^2}{\P_*}
    \right\} \ .
$$
\end{theorem}

The $h_2^{\pm}$ type bounds are the most appealing of the $p$-gradient bounds 
because the upper and lower bounds are the closest.
For instance, if $C = h_1^+(A)\,h_{\infty}^+(A) / h_2^+(A)^2$ then the gap between the upper and lower
bounds for $h_1^+$ or $h_{\infty}^+$ is at least $C$ 
(since $\frac 12\,h_1^+\,h_{\infty}^+ \geq \frac 12\,h_2^{+2}\geq \psi^+(A)$ already gives $C$
in the first inequality),
whereas the gap between the upper and lower bounds in terms of $h_2^{+2}$
is at most $\log(12\,C)$, typically a much smaller quantity. Moreover, the upper bound in terms of 
$h_2^+(A)^2$ is tighter than the upper bound for any $p\neq 2$, as can be proven via Cauchy-Schwartz.

The lower bounds on $\psi^{\pm}(A)$ for $p\neq 2$ can be considered as approximations of $\frac 12\,h_2^{+2}$,
or a bit more loosely, as approximations of $\frac 12\,h_p^+\,h_q^+$. The Jerrum-Sinclair type bound
$\psi^+(1/2)\geq \frac 12\,\P_*^{-1}\,h_1^{+2}(1/2)$ is the natural approximation to 
$\frac 12\,h_2^{+2}$ in terms of $h_1^+$, while the Alon type bound 
$\psi^+(1/2)\geq \frac 12\,\P_{min}\,h_{\infty}^{+2}(1/2)$ is natural for $h_{\infty}^+$. 
It is too much to expect the upper and lower bounds to match, 
so the extra $\log$ term in the case of $p=2$ is not much of a penalty.

Let us now look at the sharpness of Theorem \ref{thm:bounds}. 

\begin{example}
Consider the natural lazy Markov chain on the complete graph $K_n$ 
given by choosing among vertices with probability
$\frac 1{2n}$ each and holding with probability $\frac 12$ (so 
$\P(x,x)=\frac 12\,(1+\frac 1n)$ and for $y\neq x:\,\P(x,y)=\frac{1}{2n}$). 
If $A\subset K$ with size $x=\pi(A)\leq 1/2$ it follows that 
$\forall x\in A:\,\P(x,A^c)=\frac{1-x}{2}$, and therefore
$h_2^+(x)=\sqrt{(1-x)/2}$ and
$$
\frac{1-x}{4} \geq \psi^+(x) \geq \frac{1-x}{4\,\log (24/(1-x))} \geq \frac{1-x}{16}\ ,
$$
where we have used that $h_1^+(A) h_{\infty}^+(A)\leq 1$. The correct answer is 
$\psi^+(x) = \int_0^{x} \frac{t\,\frac{1-x}{2}}{x^2}\,dt = (1-x)/4$.

Likewise, $\forall x\in A^c:\,\P(x,A)=\frac{x}{2}$ and so
$h_2^-(x)=(1-x)/\sqrt{2x}$, which combined with $\psi_{gl}(A)=\psi^+(A)+\psi^-(A)$
and the bounds in Theorem \ref{thm:bounds} implies that
$$
\frac{1-x}{4x} \geq \psi_{gl}(x)
  \geq \frac{1-x}{4\,\log (24/(1-x))} + \frac{(1-x)^2}{4\,x\,\log (24/(1-x)^2)}  
  \geq \frac{1-x}{17\,x}\,,
$$ 
while the correct answer is 
$$
\psi_{gl}(x)=\int_0^{x} \frac{t\,\frac{1-x}{2}}{x^2}\,dt 
    + \int_{x}^1 \frac{(1-t)\,\frac{x}{2}}{x^2}\,dt 
 = \frac{1-x}{4x}\,.
$$

Theorem \ref{thm:main} with $\psi^+=1/32$ (as found above with $h_2^+$) implies mixing in 
$\chi^2(1/4) \leq 64\,(\log n + \log 4)$,
while Theorem \ref{thm:spectral} leads to a spectral bound of $\lambda \geq 1/128$,
which are correct orders for $\chi^2$ and $\lambda$. 
With the $\psi_{gl}$ lower bound found above we also have the correct $\tau(1/4)=O(1)$.
\end{example}

This example shows that for Markov chains with very high expansion
the bounds on $\psi_{gl}(A)$ given by Theorem \ref{thm:bounds} 
can lead to very good mixing time bounds. However, few Markov chains have such high expansion,
and so in future examples we deal only with $h_2^+(x)$ and $\psi^+(A)$.

The sharpness of
  the lower bound depends on the sharpness of the 
  $\log(h_1^+(A)\,h_{\infty}^+(A)/h_2^+(A)^2)$ term in the denominator. We give here an
  example in which the lower bound is tight, up to a factor of $1.6$, for every ratio 
  $h_2^+(A)^2/h_1^+(A)\,h_{\infty}^+(A)$ and every set size $\pi(A)$.

\begin{example}
Let the state space $K=[0,1]$ and fix some $\epsilon \leq 1/2$. For ease of
computation we consider this continuous space, but the results of this example apply to
finite spaces as well by dividing $K=[0,1]$ into states (intervals) of size $1/n$
and taking $n\rightarrow\infty$.

If $t\in[0,1/2]$ then consider the reversible Markov chain with uniform stationary
distribution on $[0,1]$ given by the transition densities
$$
\frac{\P(t,dy)}{dy} = \frac{\P(y,dt)}{dt}
= \begin{cases}
1 & if\ t \leq \epsilon \\
(\epsilon/t)^2 & if\ \epsilon < t \leq 1/2 \\
\end{cases}
\textrm{ when }
t\in\left(0,\frac 12\right) \textrm{ and } y\in\left(\frac 12,1\right),
$$
holding with the remaining probability.
Then, when $A = [0,x]\subseteq[0,1/2]$ it follows that
$$
\Psi(t,A^c) = \int_{x-t}^x \P(y,A^c)\,dy
=\begin{cases}
\frac{\epsilon^2\,t}{2\,x(x-t)} & if\ t\leq x-\epsilon \\
\frac{\epsilon(x-\epsilon)}{2x}+\frac{\epsilon-(x-t)}{2}
    & if\ t \in(x-\epsilon,\,x]
\end{cases}\,.
$$
Some computation shows that
$\psi^+(A)=\frac{\epsilon^2}{2\,x^2}\,(\frac 12+\log(x/\epsilon))$,
$h_1^+(A)=\frac{\epsilon}{2\,x}(2-\epsilon/x)$, 
$h_2^+(A)=\frac{\epsilon}{\sqrt{2}\,x}(1+\log(x/\epsilon))$
and $h_{\infty}^+(A)=1$. This leads to the relation 
\begin{equation} \label{eqn:lower}
\psi^+(A) = h_2^+(A)^2\,
     \frac{\frac 12 + \log(x/\epsilon)}
          {\left(1 + \log(x/\epsilon)\right)^2}
  < \frac{h_2^+(A)^2}{1 + \log(x/\epsilon)}\ .
\end{equation}

Let $\beta\in(0,1]$ be such that 
$$
\frac{\epsilon}{x}\,\beta^{-2} := \frac{h_2^+(A)^2}{h_1^+(A)} 
   = \frac{\epsilon}{x}\,\frac{(1+\log(x/\epsilon))^2}{2-\epsilon/x}\,.
$$
Solving for $\epsilon$ shows that $\epsilon\in(0,x]$ is the unique solution to
$
\beta= \frac{\sqrt{2-\epsilon/x}}{1-\log(\epsilon/x)}\ .
$
For every $\beta$ this gives a Markov chain and an $\epsilon$ such that 
$\epsilon/x = \beta^2\,h_2^+(A)^2/h_1^+(A)$. Letting $\beta=e^{-k/2}$ for an arbitrary $k$, and
substituting this $\epsilon$ into (\ref{eqn:lower}) gives an example with
$$
\psi^+(A) < \frac{h_2^+(A)^2}{1 + \log(x/\epsilon)}
  = \frac{h_2^+(A)^2}{k + \log\left(\frac{h_1^+(A)\,h_{\infty}^+(A)}{h_2^+(A)^2}\right)} \ .
$$
This shows that no constant $k$ in the denominator will suffice to replace the $1/2$
in Theorem \ref{thm:bounds}. 

Moreover, for every fixed $x$, as $\beta$ (or equivalently, as $\epsilon$)
varies the ratio
$\frac{h_2^+(A)^2}{h_1^+(A)\,h_{\infty}^+(A)} = \frac{\epsilon}{x}\,\beta^{-2}
  = \frac{\epsilon}{x}\,\frac{(1-\log(\epsilon/x))^2}{2-\epsilon/x}$
varies over the complete range $(0,1]$.
Therefore, for all set sizes $x$ and all ratios $h_2^{+2}/h_1^+\,h_{\infty}^+$ the
lower bound in the theorem is within a factor $2\log(12)\approx 5$ of optimal.
In fact, if the $\alpha$ quantity in the proof (see below) is optimized then
the $\alpha$ form is within a factor $1.6$ of optimal.
\end{example}

\begin{proof}[Proof of Theorem \ref{thm:bounds}]
We give the proof for the $\psi^+(A)$ case. The $\psi^-(A)$ case is similar.

First let us simplify the terms in the theorem.

Without loss, assume the state space is $K = [0,1]$ and $A=[0,x]$ where $x:=\pi(A)$. 
Order the points in $A$ in decreasing exit 
probability, so that $y,z\in A : y < z \implies \P(y,A^c) \geq \P(z,A^c)$. Then
$$
\forall t\in[0,x] : \Psi(t,A^c) = \int_{x-t}^x \P(y,A^c)\,dy \ .
$$
It follows that
$$
\psi^+(A)\,\pi(A)^2 = \int_0^x \Psi(t,A^c)\,dt = \int_0^x \int_{x-t}^x \P(y,A^c)\,dy\,dt = \int_0^x y\,\P(y,A^c)\,dy\ ,
$$
where the last equality comes from changing the order of integration.

Observe also that
$$
\Q_2(A,A^c) = \sum_{v\in A} \sqrt{\P(v,\,A^c)}\,\pi(v) = \int_0^x \sqrt{\P(t,A^c)}\,dt\ .
$$

We begin with the upper bound on $\psi^+(A)$.
\begin{eqnarray*}
\psi^+(A)\,\pi(A)^2 &=& \int_0^x \sqrt{\P(t,A^c)}\,t\,\sqrt{\P(t,A^c)}\,dt \\
  &\leq& \int_0^x \sqrt{\P(t,A^c)}\,\int_0^t \sqrt{\P(y,A^c)}\,dy\,dt \\
  &=&   \frac 12\,\left( \int_0^x \sqrt{\P(t,A^c)}\,dt\right)^2
\end{eqnarray*}
where the inequality is due to $\P(t,A^c)$ being non-increasing,
and the final equality follows from
$$
\int_0^x\,\int_0^t \sqrt{\P(t,A^c)}\,\sqrt{\P(y,A^c)} dy\,dt 
 = \int_0^x\,\int_t^x \sqrt{\P(t,A^c)}\,\sqrt{\P(y,A^c)} dt\,dy
$$
by changing the order of integration.

Now the first lower bound on $\psi^+(A)$. For any $\epsilon\in[0,x]$ 
\begin{eqnarray*}
\Q_2(A,A^c) &=& \int_0^{\epsilon} 1\cdot\sqrt{\P(t,A^c)}\,dt 
            + \int_{\epsilon}^{\Q_{\infty}(A,A^c)} \sqrt{t\,\P(t,A^c)}/\sqrt{t}\,dt \\
 &\leq& \sqrt{\epsilon\,\Q_1(A,A^c)}  + \sqrt{ \int_{\epsilon}^{\Q_{\infty}(A,A^c)} t\,\P(t,A^c)\,dt \,\int_{\epsilon}^{\Q_{\infty}(A,A^c)} \frac{dt}{t}} \\
 &\leq& \sqrt{\epsilon\,\Q_1(A,A^c)} + \sqrt{\psi^+(A)\,\pi(A)^2\,\log(\Q_{\infty}(A,A^c)/\epsilon)} 
\end{eqnarray*}
where the first inequality is by Cauchy-Schwartz.  It follows that
$$
\psi^+(A) \geq \frac{\left(\Q_2(A,A^c) - \sqrt{\epsilon\,\Q_1(A,A^c)} \right)^2}{\pi(A)^2\,\log(\Q_{\infty}(A,A^c)/\epsilon)}\ .
$$
Letting $\epsilon = \alpha^2\,\Q_2(A,A^c)^2/\Q_1(A,A^c)$ completes the lower bound. The bound
stated in the theorem follows by setting $\alpha = 1-1/\sqrt{2}$. This is a refinement of an argument of
Morris and Peres \cite{MP03.1}.

The second lower bound will be worked out for the case of general $p$-gradient
$h_p^+(A)$, so that Theorem \ref{thm:extensions} follows easily as well.
It follows from the definitions that if $B\subseteq A$ then
$\Q(A\setminus B,A^c) \geq \P_{min}^{1-1/p}\,\Q_p(A\setminus B,A^c)
\geq \P_{min}^{1-1/p}\,(\Q_p(A,A^c)-\pi(B)\,\sqrt[p]{\P_*})$. 
Then
\begin{eqnarray*}
\psi^+(A)\,\pi(A)^2 &=& \int_0^x \Q([x-t,x],A^c)\,dt \\
  &\geq& \P_{min}^{1-1/p}\,\int_0^x \max\{0,\,\Q_p(A,A^c)-t\,\sqrt[p]{\P_*}\}\,dt \\
  &=& \frac 12\,\frac{\P_{min}^{1-1/p}}{\sqrt[p]{\P_*}}\,\Q_p(A,A^c)^2\ .
\end{eqnarray*}
\end{proof}

\begin{proof}[Proof of Theorem \ref{thm:extensions}]
The upper bound for $\psi^+(A)$ follows by applying Cauchy-Schwartz to show
$h_2^+(A)^2\leq h_1^+(A)\,h_{\infty}^+(A)$, and then substituting this into
Theorem \ref{thm:bounds}. The lower bound follows from choosing the appropriate
$p$-gradients in the final section of the proof of Theorem \ref{thm:bounds}.
The bounds for $\psi^-(A)$ are proven similarly.
\end{proof}


\section{A Grid Walk} \label{sec:applications}

The previous two sections have shown that the discrete gradients provide a nice extension
of past isoperimetry results, and that the $h_2^{\pm}$  bounds 
provide relatively tight upper and lower bounds on $\psi_{gl}(A)$ and $\psi^+(A)$. 
In this section we provide an application, a near-optimal result on a random walk on 
the binary cube $2^n$, and more generally on the grid $[k]^n$. 

The quantity $h_2^+(x)$ has been studied in the theory of concentration of measure.
In particular, Talagrand \cite{Tal93.1} has shown that
$$
h_2^+(x) \geq \frac{1}{4} \sqrt{\frac{-\log x(1-x)}{n} }\ .
$$

From this and Theorem \ref{thm:bounds} we obtain the following result.

\begin{corollary} \label{cor:hypercube}
The mixing time of the lazy random walk on the cube $\{0,1\}^n$ is
$$
\tau = O\left(n\,\log^2 n\right)\ .
$$
\end{corollary}

This method (also used by Montenegro \cite{Mon02.1} and Morris and Peres \cite{MP03.1})
gives the first isoperimetric proof that $\tau$ is quasi-linear for this particular random walk.
In fact, the bound can be improved to $\tau = O(n\,\log n)$ because the proofs of both bounds in 
Theorem \ref{thm:mixingtime} only require that some expanding sequence of sets be considered.
When $\pp^{(0)}$ is a point (the worst case) then
these sets are ``fractional hamming balls'' in the blocking conductance case,
and hamming balls in the evolving set case. In both cases a modified
quantity of the form
$$
\hat{\psi}(A) = \Omega\left(\frac{-\log \pi(A)\pi(A^c)}{n}\right)
$$
can be found \cite{Mon02.1,MP03.1}. Unfortunately neither method extends to even something
as similar as the grid $[k]^n$, as the level sets are no longer Hamming balls.

We now give a proof extending Talagrand's result to the case of
the grid $[k]^n$. We will require the isoperimetric quantity
$$
\beta^{+2} = \inf_{A\subset K} \frac{\Q(A,A^c)^2}{\pi(A)\pi(A^c)}
   = \inf_{A\subset K} \Phi(A)^2\,\frac{\pi(A)}{\pi(A^c)}
$$
studied by Murali \cite{Mur01.1}.

\begin{theorem} \label{thm:tensorization}
Suppose that $K^n=K_1 \times K_2 \times \cdots K_n$ is a Cartesian product of Markov chains.
Then
$$
h_2^+(x) \geq  \frac 12\,\min \beta^+(K_i)\,\sqrt{\frac{- \log x(1-x)}{n}}\ .
$$
\end{theorem}

Houdr\'e and Tetali \cite{HT03.1} studied $h_2^+$ on products and found
that $h_2^+(K^n) \geq \frac{1}{2\sqrt{6\,n}}\,\min h_2^+(K_i)$. The advantage of
our theorem is that it considers set sizes as well, which is important for
studying mixing time.

\begin{proof}
The proof will show a chain of inequalities relating isoperimetric quantities introduced by
various authors. Our main interest is to lower bound $h_2^+(x)$, so 
rather than go into details of these quantities we will simply give definitions
and state the inequalities we need. The interested reader can learn more about these
quantities in \cite{Bob97.1,BG.1,Mur01.1}.

Bobkov's constant $b_p^+$ is defined to be the largest constant such that
for all $f : X\rightarrow[0,1]$
$$
I_{\gamma}(E f) \leq E \sqrt{I_{\gamma}^2(f) + (D_p^+ f)^2 / b_p^{+2}}\ ,
$$
where $I_{\gamma}(x)$ is the so-called Gaussian isoperimetric function,
$$
I_{\gamma}(x) = \varphi\circ \xi^{-1}(x) \quad \textrm{where} \quad \varphi(t) = \frac{1}{\sqrt{2\pi}}\,e^{-t^2/2}
 \quad \textrm{and} \quad \xi(t) = \frac{1}{\sqrt{2\pi}}\,\int_{-\infty}^t e^{-y^2/2}\,dy\ ,
$$
which makes its appearance in many isoperimetric results, and where 
$$
D_p^+ f = \left[\int \left((f(x)-f(y))^+\right)^p\,\P(x,y)\,dy\right]^{1/p} \ .
$$

Note that $I_{\gamma}(0) = I_{\gamma}(1) = 0$, $I_{\gamma}(x) \geq x(1-x)\sqrt{\log(1/x(1-x))}$,
and if $A\subset K$ then $E \left(D_p^+ 1_A\right) = \Q_p(A,A^c)$. The test 
functions $f=1_A$ thus show that
$$
b_p^+ \leq \frac{\Q_p(A,A^c)}{I_{\gamma}(\pi(A)) }
      \leq \frac{2\,h_p^+(A)}{\sqrt{-\log \pi(A)\pi(A^c)}}\,,
$$
so in a sense Bobkov's constant is a functional form of 
$\inf_{A\subset K} h_p^+(A)/\sqrt{-\log \pi(A)\pi(A^c)}$,
just as the spectral gap $\lambda$ can be considered a functional form of the
conductance $\Phi$.

It is known that Bobkov's constant $b_2^+$ tensorizes as
$$
b_2^+(K^n) = \frac{1}{\sqrt n}\,\min b_2^+(K_i)\ .
$$

Lower bounding $b_2^+(K_i)$ is difficult. One method is to use
$b_2^+(K_i) \geq b_1^+(K_i)$, which follows easily from the
definitions. In her Ph.D. Dissertation Murali \cite{Mur01.1} showed that
$b_1^+(K_i) \geq \beta^+(K_i)$. We then have the following chain of inequalities:
$$
\frac{2\,h_2^+(x)}{\sqrt{-\log x(1-x)}}
  \geq b_2^+(K^n) = \frac{1}{\sqrt n}\,\min b_2^+(K_i) \geq \frac{1}{\sqrt n}\,\min b_1^+(K_i) \geq \frac{1}{\sqrt n}\,\min \beta^+(K_i)\ .
$$
\end{proof}

Theorem \ref{thm:tensorization} is unlikely to prove of much use for Markov chains, other than ones with
relatively small state spaces, because $\beta^{+2} = \inf_{A\subset K} \Phi(A)^2\,\frac{\pi(A)}{\pi(A^c)} \leq 2\pi_0$ will be extremely
small unless $\pi_0$ is not too small. However, if a better method is found to lower bound
Bobkov's constant $b_1^+(K)$ or $b_2^+(K)$ then the theorem,
with $\beta^+(K_i)$ replaced by $b_1^+(K_i)$ or $b_2^+(K)$, could prove useful for general tensorization results.
Nevertheless, the $\beta^+$ method suffices for our purposes here.

\begin{corollary}
The mixing time of the lazy random walk on the grid $[k]^n$ is
$$
\tau = O\left(k^2\,n\,\log^2 n\right)\ .
$$
\end{corollary}

\begin{proof}
The lazy random walk on the line $[k]$ satisfies
$\Phi(x) \geq 1/(4k\,x)$. It follows that $\beta^{+2} \geq 1/4k^2$.
By Theorem \ref{thm:tensorization} $h_2^+(x)\geq \frac{1}{4\,k\sqrt{n}}\sqrt{-\log x(1-x)}$.
Combining this with Theorems \ref{thm:mixingtime} and \ref{thm:bounds},
and the simplification $\frac{h_1^+(A)\,h_{\infty}^+(A)}{h_2^+(A)^2}\leq \frac{\P_*}{\P_{min}}$
discussed after Theorem \ref{thm:bounds} gives $\tau=O(k^2\,n\,\log n\,\log \log k^n)$.

The $\log \log k^n$ term can be improved to $\log n$ because of
``ultracontractivity'' of geometric Markov chains, 
such as this one on $[k]^n$, which implies that $\pi_0$ may be taken as $2^{-n}$ 
instead of $k^{-n}$ with only the addition of a small constant factor. \cite{Mon03.1}
\end{proof}

This is the first isoperimetric proof of $\tau = O^*(k^2\,n)$ for this Markov chain,
and is quite close to the correct $\tau = \Theta(k^2\,n\,\log n)$.


\section*{Acknowledgments}

The author thanks Gil Kalai for introducing him to Talagrand's result,
which led to the author's study of discrete gradients.
Christian Houdr\'e, Ravi Kannan and L\'aszl\'o Lov\'asz were of
great help at various stages of this work. Suggestions of the
referees helped to greatly improve the clarity of this paper.


\bibliographystyle{plain}
\bibliography{../references}

\end{document}